\documentclass[11pt]{article}
\usepackage{mathrsfs}
\usepackage{amsmath}
\usepackage{amsfonts}
\usepackage{amsfonts,amssymb}
\usepackage{amscd}
\usepackage{epsf}
\usepackage{amsmath,amssymb,latexsym,color}
\newtheorem{thm}{Theorem}[section]
\newtheorem{prop}[thm]{Proposition}
\newtheorem{lem}[thm]{Lemma}

\newtheorem{example}{Example}[section]

\newtheorem{remark}{Remark}[section]

\newcommand{\proof}{{\it Proof.\quad}}
\newcommand{\qed}{\hfill\Box\medskip}

\usepackage{CJK}
\begin{document}

\title{\bf Quasi Regular Semilattice and Association Schemes in Singular Linear Space }

\author{ Zhang Baohuan$^1$\quad Yue Mengtian$^2$
  \quad Li Zengti$^1$\thanks{lizengti@126.com} \\
{\footnotesize  \em 1.Math. and  Inf. College, Langfang Teachers'
College, Langfang, 065000,  China }\\
{\footnotesize  \em
2.Department of Scientific Research, Langfang Teachers'
College, Langfang, 065000,  China }}

\date{}
 \maketitle


{\noindent\small{\bf Abstract}\quad Let $\mathbb{F}_q^{n+l}$ denote the $(n+l)$-dimensional singular linear space over a finite field $\mathbb{F}_q$. For a fixed integer $m\leq\min\{n,l\}$, denote by $\mathcal{L}^m_o(\mathbb{F}_q^{n+l})$ the set of all subspaces of type $(t,t_1)$, where $t_1\leq t\leq m$.
Partially ordered by ordinary inclusion, one family of quasi regular semilattices is obtained. Moreover, we obtain a association schemes and discuss the bound of a $M$-clique.



\section{Introduction}
It is well known that lattice is an important part of poset's theory. its theory play an important role in many branches of mathematics, such as
computer logical design.
The results on the
lattices generated by transitive sets of subspaces under finite
classical groups may be found in Huo, Liu and
Wan \cite{Huo1,Huo2,Huo3}. In \cite{guo}, Guo discussed the lattices
associated with finite vector spaces and finite affine spaces. In \cite{DD}, P. Delsarte discussed the regular semilattices
in finite vector spaces. In this paper, we obtain a new quasi regular semilattice and a new association scheme in the singular linear space. More over we discuss the bounds of a $M$-clique.

The rest of this paper is organized as followed. In section 2, we discuss some definitions and terminologies about  lattices, regular semilattices, association schemes and $M$-cliques. In section 3, we construct a family of quasi regular semilattice, and then compute its parameters. In section 4, we obtain a new association. In section 5, we discuss the bound of a $M$-clique.

\section{Preliminaries}

Let $(P,\leq)$ be a poset. We write $a<b$ whenever $a\leq b$ and
$a\neq b$.  If $P$ has the minimum (respectively
maximum) element, then we denote it by 0 (respectively $\upharpoonleft$), and
say that $P$ is a poset with 0 (respectively $\upharpoonleft$).  A poset $P$ is said to be a  {\it semilattice} if $a\wedge b:=\rm{inf}\{a,b\}$ exist for any two elements $a, b\in P$. Let P be a finite poset with 0. If there is a function $r$ from $P$ to set
of all the nonnegative integers such that
\begin{itemize}
  \item [\rm(1)] r(0)=0,
  \item [\rm(2)] $r(b)=r(a)+1$, if $a\lessdot b$.
\end{itemize}
Then $r$ is said to be the  {\it rank function} on $P$.
Note that the rank function on $P$ is unique if it exists.

Let $P$ be a semilattice, and let $P=X_0\cup X_1\cup\cdots\cup X_m$, where $X_i =\{x\in P| r(x) = i\}, i = 0, 1,\cdots, m.$
 The semilattice $(P, \leq)$ is called {\it regular} if  the following three properties hold:
\begin{itemize}
  \item [\rm(i)] Given $y \in X_m , z \in X_r$ with $z \leq y$, the number of points $u \in X_s$
such that $z\leq u \leq y$ is a constant $\mu(r, s)$.
  \item [\rm(ii)] Given $u \in X_s$ , the number of points $z \in X_r$, such that $z \leq u$ is a
constant $\nu(r,s)$.
  \item [\rm(iii)] Given $a\in X_r, , y \in X_m$, with $a\wedge y \in X_j$, the number of pairs
$(b, z) \in X_s\times X_m$ such that $b \leq z, b \leq y, a \leq z$ is a constant $\pi(j, r, s).$
\end{itemize}
In this paper, we define the concept of quasi regular semilattice as follows.

The semilattice $(P, \leq)$ is called {\it quasi regular} if  the following three properties hold:
\begin{itemize}
  \item [\rm(i)] Given $y \in X^{m'}_m , z \in X^{r'}_r$ with $z \leq y$, the number of points $u \in X^{s'}_s$
such that $z\leq u \leq y$ is a constant $\mu(r(r'), s(s');m')$.
  \item [\rm(ii)] Given $u \in X^{s'}_s$ , the number of points $z \in X^{r'}_r$, such that $z \leq u$ is a
constant $\nu(r(r'),s(s'))$.
  \item [\rm(iii)] Given $a\in X^{r'}_r, y \in X^{m_1}_m$, with $a\wedge y \in X^{j'}_j$, the number of pairs
$(b, z) \in X^{s'}_s\times X^{m'}_m$ such that $b \leq z, b \leq y, a \leq z$ is a constant $\pi(j(j'), r(r'), s(s');m').$
\end{itemize}
Here $X_i=X_i^{0}\cup X_i^{1}\cup\cdots\cup X_i^{i}$, and $X_i^{j}\cap X_i^{k}=\emptyset$ for $j\neq k.$

Let $X$ be a finite set of vertices. A $d$-class association scheme on $X$ consists of a set of $d+1$
symmetric relations $R_0,R_1,\cdots,R_d$ on $V$, with identity relation $R_0=\{(x,x)|x\in V\},$
such that any two vertices are in precisely one relation. Denoted by $(X,\{R_i\}_{0\leq i\leq d})$. Furthermore, there are intersection
numbers $p^k_{ij}$ such that for any $(x,y)\in R_k$, the number of vertices $z$ such that $(x,z)\in R_i$ and $(z,y)\in R_j$ equals $p^k_{ij}$.

The nontrivial relations can be considered as graphs, which in our case are undirected.
One immediately sees that the respective graphs are regular with degree $v_i=p_{ii}^0$. For the
corresponding adjacency matrices $A_i$ the axioms of the scheme are equivalent to
$$\sum_{i=0}^dA_i=J,\ A_0=I,\ A_i=A_i^T,\ A_iA_j=\sum_{k=0}^dp_{ij}^kA_k.$$
It follows that the adjacency matrices generate a $(d+1)$-dimensional commutative algebra
$\mathfrak{A}$ of symmetric matrices. This algebra was first studied by Bose and Mesner \cite{BM} and is
called the Bose-Mesner algebra of the scheme.

A nonzero vector $\alpha\in \mathbb{R}^{|X|}$ is said to be {\it $\theta$-positive },  if $C_i(\alpha)=v_i^{-1}<\alpha,A_i\alpha>\geq 0$ holds for every $i$. Let $M$ be a subset of $[0,n]$, containing $0$, Any positive vector $\alpha$ is called an {\it $M$-clique} if it satisfies
$$C_j(\alpha)=v_j^{-1}<\alpha,A_j\alpha>=0,\ for\ all\ j\not\in M$$

Let $x_0$ be a point of $X$ and  $s$ be an integer in $[0,n]$, A given $\theta$-positive vector $\alpha\in \mathbb{R}^{|X|}$ is called a {\it unicoloured vector} of center $x_0$ and colour $s$, if the following condition holds:
$$\alpha(x)=0,\ unless\ (x,x_0)\in R_s$$

Let $\mathbb{F}_q$ be a finite field with $q$ elements, where $q$ is a prime power. For two non-negative integers
$n$ and $l$, $\mathbb{F}_q^{n+l}$ denotes the $(n +l)$-dimensional row vector space over $\mathbb{F}_q$. The set of all $(n +l)\times(n +l)$
nonsingular matrices over $\mathbb{F}_q$ of the form
$$\left(\begin{array}{cc}
T_{11}&T_{12}\\
0&T_{22}
\end{array}\right),$$
where $T_{11}$ and $T_{22}$ are nonsingular $n\times n$ and $l\times l$ matrices, respectively, forms a group under matrix
multiplication, called the singular general linear group of degree $n + l$ over $\mathbb{F}_q$ and denoted by
$GL_{n+l,n}(\mathbb{F}_q)$.

Let $P$ be an $m$-dimensional subspace of $\mathbb{F}_q^{n+l}$, denote also by $P$ a $m\times (n + l)$ matrix of rank $m$
whose rows span the subspace $P$ and call the matrix $P$ a matrix representation of the subspace $P$ .
There is an action of $GL_{n+l,n}(\mathbb{F}_q)$ on $\mathbb{F}_q^{n+l}$ defined as follows
$$\begin{array}{rrr}
\mathbb{F}_q^{n+l}\times GL_{n+l,n}(\mathbb{F}_q)&\longrightarrow &\mathbb{F}_q^{n+l},\ \ \ \ \ \ \ \ \ \ \ \ \ \ \ \ \ \ \ \ \ \ \ \ \ \ \ \ \ \ \ \ \ \ \ \ \ \ \ \ \\
((x_1,\cdots,x_n,x_{n+1},\cdots,x_{n+l}),T)&\longmapsto&(x_1,\cdots,x_n,x_{n+1},\cdots,x_{n+l})T.
\end{array}$$
The above action induces an action on the set of subspaces of $\mathbb{F}_q^{n+l}$; i.e., a subspace $P$ is carried by
$T\in GL_{n+l,n}(\mathbb{F}_q)$ to the subspace $PT$. The vector space $\mathbb{F}_q^{n+l}$ together with the above group action is
called the $(n +l)$-dimensional singular linear space over $\mathbb{F}_q$.

For $1\leq i \leq n + l$, let $e_i$ be the row vector in $\mathbb{F}_q^{n+l}$ whose $i$-th coordinate is $1$ and all other coordinates
are $0$. Denote by $E$ the $l$-dimensional subspace of  $\mathbb{F}_q^{n+l}$ generated by $e_{n+1}, e_{n+2}, \cdots, e_{n+l}$.
A $m$-dimensional subspace $P$ of  $\mathbb{F}_q^{n+l}$ is called a subspace of type $(m,k)$ if $\dim(P \cap E) = k$.

For a fixed subspace $P$ of type $(m_1,k_1)$ in $\mathbb{F}_q^{n+l}$, let $\mathcal{M}'(m_1,k_1;m,k;n+l,n)$ denote the set of all the
subspaces of type $(m,k)$ containing $P$.  $\mathcal{M}'(m_1,k_1;m,k;n+l,n)$ is non-empty if and only if $0\leq k_1\leq k\leq l$ and $0\leq m_1-k_1\leq m-k\leq n$. Let $N'(m_1,k_1;m,k;n+l,n)=|\mathcal{M}'(m_1,k_1;m,k;n+l,n)|.$ The formula is given by Kaishun Wang, Jun Guo, Fenggao Li~(see \cite{wgl}) as follows.
$$N'(m_1,k_1;m,k;n+l,n)=q^{(l-k)(m-k-m_1+k_1)}\left[n-(m_1-k_1)\atop (m-k)-(m_1-k_1) \right]_q\left[l-k_1\atop k-k_1\right]_q.$$

\begin{lem}~\label{l2-l}
 Let $0\leq k_1\leq k_2\leq k\leq l,0\leq m_1-k_1\leq m_2-k_2\leq m-k\leq n.$ For $i=1,2,3$, assume that $U_i$ be a set of all the subspaces of type $(m_i,k_i)$ in $\mathbb{F}_q^{n+l}$, and  $U$ be  a set of all the subspaces of type $(m,k)$ in $\mathbb{F}_q^{n+l}$ such that $U_1\subseteq U_2\subseteq U$ and $U_3\subseteq U$. if $U_3\cap U_2 = U_1$, then the number of $U_3$ is
$$q^{(\delta_3-\delta_1)(\delta_2-\delta_1+k-k_3)+(k_3-k_1)(k_2-k_1)}\left[\delta-\delta_2\atop \delta_3-\delta_1\right]_q\left[k-k_2\atop k_3-k_1\right]_q,$$
where $\delta=m-k,\delta_i=m_i-k_i,i=1,2,3$.
\end{lem}
\proof Let $\delta_1=m_1-k_1,\delta_2=m_2-k_2,\delta_3=m_3-k_3,\delta=m-k$.
By the transitivity of $GL_{n+l,n}(\mathbb{F}_q)$ on the set of subspaces of the same type, we may assume that
$$\begin{array}{l}
                              \begin{array}{llll}
                        \ \ \ \ \ \ \ \ \  \delta_1&n-\delta_1&k_1& l-k_1
                               \end{array}       \\
                              U_1=
                            \left(\begin{array}{llll}
                           I\ \ \ \ &0\ \ \  &0\ \ &0\\
                           0\ &0\ \ \ &I\ \ \ &0
                              \end{array}\right)
                       \end{array}
                        \begin{array}{l} \begin{array}{l}
                        \par \\
                             \delta_1\\
                             k_1
                       \end{array}  \end{array}
$$
$$\begin{array}{l}
                              \begin{array}{llllll}
                        \ \ \ \ \ \ \ \ \  \delta_1&\delta_2-\delta_1&n-\delta_2&k_1&k_2-k_1&l-k_2
                               \end{array}       \\
                              U_2=
                            \left(\begin{array}{llllll}
                           I\ \ \ \  &0\ \ \ \ \ \  &0\ \ \ \ \ &0\ \ \ \ \  &0\ \ \ \  &0\\
                           0\ &0&0\ \ \ &I\ \ \ &0&0\\
                        0&I&0&0&0&0\\
                        0&0&0&0&I&0
                              \end{array}\right)
                       \end{array}
                        \begin{array}{l} \begin{array}{l}
                        \par \\
                             \delta_1\\
                             k_1\\
                            \delta_2-\delta_1\\
                             k_2-k_1
                       \end{array}  \end{array}
$$
and
$$\begin{array}{l}
                              \begin{array}{llllllll}
                        \ \ \ \ \ \ \ \ \ \  \delta_1&\delta_2-\delta_1&\delta-\delta_2&n-\delta&k_1&k_2-k_1&k-k_2&l-k
                               \end{array}       \\
                              U=
                            \left(\begin{array}{llllllll}
                           I\ \ \ \ &0\ \ \ \ \ \ \ &0\ \ \ \ \  &0\ \ \ \ \ &0\ \ \ \ \ \ &0\ \ \ \ \ \  \ \ &0\ \ \ \ \  &0\\
                           0\ &0&0\ &0 &I\ \ \ &0&0&0\\
                        0&I&0&0&0&0&0&0\\
                        0&0&0&0&0&I&0&0\\
                        0&0&I&0&0&0&0&0\\
                        0&0&0&0&0&0&I&0\\
                              \end{array}\right)
                       \end{array}
                        \begin{array}{l} \begin{array}{l}
                        \par \\
                             \delta_1\\
                             k_1\\
                             \delta_2-\delta_1\\
                             k_2-k_1\\
                             \delta-\delta_2\\
                             k-k_2
                       \end{array}  \end{array}
$$
Since
$U_3\cap U_2 = U_1,$
we have
$$\begin{array}{l}
                              \begin{array}{llllllll}
                        \ \ \ \ \ \ \ \ \ \delta_1&\delta_2-\delta_1&\delta-\delta_2&n-\delta&k_1&k_2-k_1&k-k_2&l-k
                               \end{array}       \\
                              U_3=
                            \left(\begin{array}{llllllll}
                           I\ \ \ \ &0\ \ \ \ \ \ &0\ \ \ \ \ \ \ \  &0\ \  &0\ \ \ \  &0\ \ \ \ \ \ \ &0\ \ \ \ \  &0\\
                           0\ &0&0\ &0 &I\ \ \ &0&0&0\\
                        0&u_{32}&u_{33}&0&0&u_{36}&u_{37}&0\\
                        0&0&0&0&0&u_{46}&u_{47}&0\\
                         \end{array}\right)
                       \end{array}
                        \begin{array}{l} \begin{array}{l}
                        \par \\
                             \delta_1\\
                             k_1\\
                             \delta_3-\delta_1\\
                             k_3-k_1\\
 \end{array}  \end{array}
$$
where  rank $u_{47}=k_3-k_1$. Note that there are $\left[k-k_2\atop k_3-k_1\right]_q$ choices for $u_{47}$.
 By the transitivity of $GL_{n+l,n}(\mathbb{F}_q)$ on the set of subspaces of the same type, the number of $U_3$'s
 does not depend on the particular choice of $u_{47}$. Pick $u_{47}=(I^{(k_3-k_1)},0)$. Then $U_3$ has a matrix representation
$$\begin{array}{l}
                              \begin{array}{lllllllll}
                        \ \  \delta_1&\delta_2-\delta_1&\delta-\delta_2&n-\delta&k_1&k_2-k_1&k_3-k_1&k-k_2-k_3+k_1&l-k
                               \end{array}       \\

                            \left(\begin{array}{lllllllll}
                           I\ \ \ \ &0\ \ \ \ \ \ \ &0\ \ \ \ \ \ &0\ \ \  &0\ \ \ \ \ \ &0\ \ \ \ \ \ \ &0\ \ \ \ \ \ \ \ \ \ \ \ \ &0\ \ \ \ \ \ \ \ \ \ \ \ &0\\
                           0\ &0&0\ &0 &I\ \ \ &0&0&0&0\\
                        0&u_{32}&u_{33}&0&0&u_{36}&0&u'_{37}&0\\
                        0&0&0&0&0&u_{46}&I&0&0\\
                         \end{array}\right)
                       \end{array}
                        \begin{array}{l} \begin{array}{l}
                        \par \\
                             \delta_1\\
                             k_1\\
                             \delta_3-\delta_1\\
                             k_3-k_1\\
 \end{array}  \end{array}
$$
where rank $u_{33}=\delta_3-\delta_1$.
Therefore the number of $U_3$ is equal to
$$q^{(\delta_3-\delta_1)(\delta_2-\delta_1+k-k_3)+(k_3-k_1)(k_2-k_1)}\left[\delta-\delta_2\atop \delta_3-\delta_1\right]_q\left[k-k_2\atop k_3-k_1\right]_q.$$
$\qed$\\

\section{The quasi regular semilattice }
For a fixed integer $m\leq\min\{n,l\}$, denote by $\mathcal{L}^m(\mathbb{F}_q^{n+l})$ the set of all subspaces of type $(t,t_1)$, where $t_1\leq t\leq m$.
 If we
partially order $\mathcal{L}^m(\mathbb{F}_q^{n+l})$ by the ordinary inclusion, then $\mathcal{L}^m(\mathbb{F}_q^{n+l})$ is a semilattice, denoted by
$\mathcal{L}^m_o(\mathbb{F}_q^{n+l})$. For any $A \in \mathcal{L}^m_o(\mathbb{F}_q^{n+l})$, The rank function of $\mathcal{L}_o(\mathbb{F}_q^{n+l})$ is defined as follows
$$r(A) =\dim (A).$$
Let
$$X_i=\{B\in \mathcal{L}^m_o(\mathbb{F}_q^{n+l})|r(B)=i\},$$
and
$$X_i^j=\{B\in X_i|\dim(B\cap E)=j\},j=0,1,\cdots,i,$$
where $E=\langle e_{n+1},e_{n+2},\cdots,e_{n+l}\rangle\in \mathbb{F}_q^{n+l}.$
We will prove that $\mathcal{L}^m_o(\mathbb{F}_q^{n+l})$ is a quasi regular semilattice and
 compute its parameters.

\begin{lem}~\label{l3-1}
 Let $A\in X^{r_1}_r$, $C\in X^{m_1}_m$ and $A\leq C$.
Then the number of $B\in X^{s_1}_s$ such that $A\leq B\leq C$ is equal to
$$\mu(r(r_1),s(s_1);m_1)=q^{(s-s_1-r+r_1)(m_1-s_1)}
\left[\begin{array}{l}
m-r+r_1-m_1\\
 s-s_1-r+r_1
 \end{array}\right]_q
 \left[\begin{array}{l}
m_1-r_1\\
 s_1-r_1
 \end{array}\right]_q
 $$
\end{lem}
\proof
 By the transitivity of $GL_{n+l}(\mathbb{F}_q)$ on the set of subspaces of the same type, we may assume that
$$\begin{array}{l}
                              \begin{array}{llll}
                        \ \ \ \ \ \ \ \  r-r_1&n-r+r_1&r_1& l-r_1
                               \end{array}       \\
                              A=
                            \left(\begin{array}{llll}
                           I\ \ \ \ \ \ \ \ \ \ \ &0\ \ \ \ \ \ \ \ &0\ \ &0\\
                           0\ &0\ \ \ &I\ \ \ &0
                              \end{array}\right)
                       \end{array}
                        \begin{array}{l} \begin{array}{l}
                        \par \\
                             r-r_1\\
                             r_1
                       \end{array}  \end{array}
$$
and
$$\begin{array}{l}
                              \begin{array}{llllll}
                        \ \ \ \ \ \ \ \  r-r_1&m-m_1-r+r_1&n-m+m_1&r_1&m_1-r_1&l-m_1
                               \end{array}       \\
                              C=
                            \left(\begin{array}{llllll}
                           I\ \ \ \ \ \ \ \ \ \ \ \  &0\ \ \ \ \ \ \ \ \ \ \ \ \ \ \ \ \ \  \ &0\ \ \ \ \ \ \ \ \  &0\ \ \ \ \ \ &0\ \ \ \ \ &0\\
                           0\ &0&0\ \ \ &I\ \ \ &0&0\\
                        0&I&0&0&0&0\\
                        0&0&0&0&I&0
                              \end{array}\right)
                       \end{array}
                        \begin{array}{l} \begin{array}{l}
                        \par \\
                             r-r_1\\
                             r_1\\
                             m-m_1-r+r_1\\
                            m_1-r_1
                       \end{array}  .\end{array}
$$
Since
$$A\leq B\leq C,$$
we have
$$\begin{array}{l}
                              \begin{array}{llllll}
                        \ \ \ \ \ \ \ \  r-r_1&m-m_1-r+r_1&n-m+m_1&r_1&m_1-r_1&l-m_1
                               \end{array}       \\
                              B=
                            \left(\begin{array}{llllll}
                           I\ \ \ \ \ \ \ \ \ \ \ \ \ \ \ \ \ &0\ \ \ \ \ \ \ \ \ \ \ \ \ \  &0\ \ \ \ \ \ \ \ &0\ \ \ \ \ \ \ &0\ \ \ \ \ \  &0\\
                           0\ &0&0\ \ \ &I\ \ \ &0&0\\
                        0&u_{32}&0&0&u_{35}&0\\
                        0&0&0&0&u_{45}&0
                              \end{array}\right)
                       \end{array}
                        \begin{array}{l} \begin{array}{l}
                        \par \\
                             r-r_1\\
                             r_1\\
                             s-s_1-r+r_1\\
                            s_1-r_1
                       \end{array} \end{array},
$$
where rank $u_{45}=(s_1-r_1)$. Note that there are $\left[m_1-r_1\atop s_1-r_1\right]_q$ choices for $u_{45}$.
By the transitivity of $GL_{n+l}(\mathbb{F}_q)$ on the set of subspaces of the same type, the number of $B$'s does not depend on the particular choice of $u_{45}$. Pick
$u_{45}=(I^{(s_1-r_1)},0)$. Then $B$ has a matrix representation
  $$\begin{array}{l}
                              \begin{array}{lllllll}
                         \ \ \delta_2&\delta_1-\delta_2&n-\delta_1&r_1&s_1-r_1&m_1-s_1&l-m_1
                               \end{array}       \\
                        \left(\begin{array}{lllllll}
                           I\ \ \ \ \ \  &0\ \  \ \ \   &0\ \ \ \ &0\ \ \ \ \  &0\ \ \ \ \ \ \ \  &0\ \ \ \ \ \ &0\\
                           0\ &0&0\ \ \ &I\ \ \ &0&0&0\\
                        0&u_{32}&0&0&0&u'_{35}&0\\
                        0&0&0&0&I&0&0
                              \end{array}\right)
                       \end{array}
                        \begin{array}{l} \begin{array}{l}
                        \par \\
                             \delta_2\\
                             r_1\\
                             \delta_3-\delta_2\\
                            s_1-r_1
                       \end{array} . \end{array}
$$
where $\delta_1=m-m_1,\delta_2=r-r_1,\delta_3=s-s_1$.
Therefore the number of $B$ is equal to
$$q^{(s-s_1-r+r_1)(m_1-s_1)}
\left[\begin{array}{l}
m-r+r_1-m_1\\
 s-s_1-r+r_1
 \end{array}\right]_q
 \left[\begin{array}{l}
m_1-r_1\\
 s_1-r_1
 \end{array}\right]_q.
 $$
 $\qed$

 \begin{lem}~\label{l3-2}
If $B\in X^{s_1}_s$, then the number of $A\in X^{r_1}_r$ such that $A\leq B$ is equal to
 $$\nu(r(r_1),s(s_1))=q^{(r-r_1)(s_1-r_1)}\left[s-s_1\atop r-r_1\right]_q\left[s_1\atop r_1\right]_q.$$
 \end{lem}
 \proof
By the transitivity of $G_{n+l,l}(\mathbb{F}_q)$ on the set of subspaces of the same type,  we may assume that
$$\begin{array}{l}
                              \begin{array}{llll}
                        \ \ \ \ \ \ \ \  s-s_1&n-s+s_1&s_1& l-s_1
                               \end{array}       \\
                              B=
                            \left(\begin{array}{llll}
                           I\ \ \ \ \ \ \ \ \ \ &0\ \ \ \ \ \ \ \ &0\ \ &0\\
                           0\ &0\ \ \ &I\ \ \ &0
                              \end{array}\right)
                       \end{array}
                        \begin{array}{l} \begin{array}{l}
                        \par \\
                             s-s_1\\
                             s_1
                       \end{array} \end{array}.
$$
Since
$$A\leq B,$$
we have
$$\begin{array}{l}
                              \begin{array}{llll}
                        \ \ \ \ \ \ \ \  s-s_1&n-s+s_1&s_1& l-s_1
                               \end{array}       \\
                              A=
                            \left(\begin{array}{llll}
                           u_{11}\ \ \ \ \ \ \ \ &0\ \ \ \ \ &u_{13}&0\\
                           0\ &0\ \ \ &u_{23}\ \ \ &0
                              \end{array}\right)
                       \end{array}
                        \begin{array}{l} \begin{array}{l}
                        \par \\
                             r-r_1\\
                             r_1
                       \end{array}  \end{array},
$$
where rank $u_{23}=r_1$. Note that there are $\left[s_1\atop r_1\right]_q$ choices for $u_{23}$.
By the transitivity of $GL_{n+l}(\mathbb{F}_q)$ on the set of subspaces of the same type, the number of $A$'s does not depend on the particular choice of $u_{23}$. Pick
$u_{23}=(I^{(r_1)}\ 0)$. Then $A$ has a matrix representation
$$\begin{array}{l}
                              \begin{array}{lllll}
                         s-s_1&n-s+s_1&r_1&s_1-r_1& l-s_1
                               \end{array}       \\
                         \left(\begin{array}{lllll}
                           u_{11}\ \ \ \ \ \ \ \ &0\ \ \ \ \ \ &0\ \ \ \ &u'_{13}\ \ \ \ \ &0\\
                           0\ &0\ \ \ &I&0\ \ \ &0
                              \end{array}\right)
                       \end{array}
                        \begin{array}{l} \begin{array}{l}
                        \par \\
                             r-r_1\\
                             r_1
                       \end{array}  \end{array}.
$$
Therefore the number of subspace $A$ is equal to
$$\nu(r(r_1),s(s_1))=q^{(r-r_1)(s_1-r_1)}\left[s-s_1\atop r-r_1\right]_q\left[s_1\atop r_1\right]_q.$$
$\qed$

\begin{lem}~\label{l3-3}
Let $A\in X^{r_1}_r$ and $B\in X^{m_2}_m$. Assume that $A\wedge B\in X^{j_1}_j$, and $(C,D)\in X_s^{s_1}\times X_m^{m_1}$
. If
$C\leq D,C\leq B,A\leq D$, then the number of $(C,D)$ is equal to
$$\begin{array}{ll}
&\pi(j(j_1),r(r_1),s(s_1);m_1)\\
=&\sum_{0\leq i\leq j,0\leq i_1\leq \min\{i,j_1\}}q^{(s-s_1-i+i_1)(j-j_1-i+i_1+m_1-s_1)+(s_1-i_1)(j_1-i_1)}\\
&\times
\left[(m-m_1)-(j-j_1)\atop (s-s_1)-(i-i_1)\right]_q\left[m_1-j_1\atop s_1-i_1\right]_qN'(r+s-i,r_1+s_1-i_1;m,m_1;n+l,n).
\end{array}
$$
\end{lem}
\proof
Since  $A\leq D,C\leq B$, we have
$$A\wedge C\leq A\wedge B\leq D.$$
Since
$$C\leq D,$$
and
$$C\wedge(A\wedge B)=A\wedge(B\wedge C)=A\wedge C\in X^{i_1}_i,0\leq i\leq j, 0\leq i_1\leq \min\{i,j_1\}.$$
 For $i\in [0,j]$ and $i_1\in [0,\min\{i,j_1\}]$, by Lemma~\ref{l2-l}, the number of $C$ is equal to
$$q^{(s-s_1-i+i_1)(j-j_1-i+i_1+m_1-s_1)+(s_1-i_1)(j_1-i_1)}
\left[(m-m_1)-(j-j_1)\atop (s-s_1)-(i-i_1)\right]_q\left[m_1-j_1\atop s_1-i_1\right]_q.$$

Since $C\leq D,A\leq D$, $C+A\leq D$. It follows from  $A\wedge C\in X_i^{i_1}$ that
 $A+C$ is a subspace of type $(r+s-i,r_1+s_1-i_1)$. Therefore, the number of $D$ is equal to $N'(r+s-i,r_1+s_1-i_1;m,m_1;n+l,n)$.
Hence the desired result follows.
 $\qed$

By Lemma~\ref{l3-1},Lemma~\ref{l3-2} and Lemma~\ref{l3-3}, we obtain the following theorem.
 \begin{thm}~\label{t3-1}
Semilattice $\mathcal{L}_0^m(\mathbb{F}_q^{n+l})$ is a quasi regular semilattice. Its parameters are given by the formulas
$$\mu(r(r_1),s(s_1);m_1)=q^{(s-s_1-r+r_1)(m_1-s_1)}
\left[\begin{array}{l}
m-r+r_1-m_1\\
 s-s_1-r+r_1
 \end{array}\right]_q
 \left[\begin{array}{l}
m_1-r_1\\
 s_1-r_1
 \end{array}\right]_q,$$
 $$\nu(r(r_1),s(s_1))=q^{(r-r_1)(s_1-r_1)}\left[s-s_1\atop r-r_1\right]_q\left[s_1\atop r_1\right]_q,$$
 and
 $$\begin{array}{ll}
&\pi(j(j_1),r(r_1),s(s_1);m_1)\\
=&\sum_{0\leq i\leq j,0\leq i_1\leq i}q^{(s-s_1-i+i_1)(j-j_1-i+i_1+m_1-s_1)+(s_1-i_1)(j_1-i_1)}\\
&\times
\left[(m-m_1)-(j-j_1)\atop (s-s_1)-(i-i_1)\right]_q\left[m_1-j_1\atop s_1-i_1\right]_qN'(r+s-i,r_1+s_1-i_1;m,m_1;n+l,n).
\end{array}
$$
\end{thm}

\section{Association schemes}

Given $d$ and $d_0$ with $0\leq d_0\leq d\leq \min\{n,l\}$. For $h=\min\{k,d_0\},$ we define
$$\widetilde{X_k}=\left\{
\begin{array}{ll}
X_k^0\cup \cdots \cup X_k^h,&if\  0\leq k\leq d-1,\\
X_d^{d_0},&if\  k=d,
\end{array}
\right.$$
and
$$R_{k} = \{(B,C)\in X_d^{d_0}\times X_d^{d_0}|B\wedge C\in \widetilde{X_{d-k}}\}.$$
Then the $R=\{R_0,R_1,\cdots,R_d\}$ is a partition of $X_d^{d_0}\times X_d^{d_0}$.
For any $B,C\in X_d^{d_0}$,
let
$$A_{k(k')}(B,C)=\left\{
\begin{array}{ll}
1,&if\ B\wedge C\in X^{k'}_{d-k},\\
0,& otherwise.
\end{array}
\right.$$
Then, for $h=\min\{d-k,d_0\}$, we define
 $$A_{k}=\left\{
\begin{array}{ll}
A_{k(0)}+A_{k(1)}+\cdots+A_{k(h)},&if\  1\leq k\leq d,\\
A_{0(d_0)},&if\  k=0.
\end{array}
\right.$$
Let
$$D_{i,k}(B,C)=\left\{
\begin{array}{ll}
1,&if\  B\leq C,\\
0,& otherwise.
\end{array}
\right.
$$
The $D_{i,k}$ is a Riemann matrix indexed by $\widetilde{X_i}\times \widetilde{X_k}$.
For convenience, we write $D_i=D_{i,d}$.
Let $\mathfrak{A}=\langle A_0,A_2,\cdots,A_d\rangle$ denote the $(d+1)$-dimensional real vector space generated by the $A_0,A_2,\cdots,A_d$.
Let $C_i=D^T_iD_i$, for  $i=0,1,\cdots,d.$
We have the lemma as follows.
\begin{lem}~\label{4l-1}
The matrices $C_0,C_1,\cdots,C_d$ generate  $\mathfrak{A}$, with
\begin{equation}\label{el1}
    C_t=\sum_{h=0}^{\min\{t,d_0\}}\nu(t(h),d(d_0))A_{0}+
    \sum_{k=1}^{d}\sum_{\lambda=0}^{\min\{d-k,d_0\}}\sum_{h=0}^{\min\{t,d_0\}}\nu(t(h),(d-k)(\lambda))A_{k},
\end{equation}
where $ t=0,1,\cdots,d.$
\end{lem}
\proof
For any $(A,B)\in X_d^{d_0}\times X_d^{d_0}$, the number of $F\in \widetilde{X_t}$ such that $F\leq A\wedge B$ is the $(A,B)$-entry in two members of (\ref{el1}).
It implies that the system (\ref{el1}) holds. On the other hand, its matrix has rank $d+1$. Hence, the result is obtained.
$\qed$
\begin{lem}~\label{4l-2}
Let $A\in X^{d_1}_d,B\in X^{s_1}_s$, then the number of $D\in X^{r_1}_r$ such that $D\leq A$ and $B\wedge D\in X^{j_1}_j$ is a constant  $\psi(j(j_1),r(r_1),s(s_1);d_1)$.
\end{lem}
\proof
For a given $k$ and $k_1$, with $0\leq k\leq h=\min\{r,s\}$ and $0\leq k_1\leq k$. We assume that $(C,D)\in (X_k^{k_1},X_r^{r_1})$ and $C\leq B, C\leq D\leq A.$
 Counting the number of the pairs $(C,D)$ in two different ways, we obtain
\begin{equation}\label{el2}
    \sum_{j=0}^h\sum_{j_1=0}^j\nu(k(k_1),j(j_1))\psi(j(j_1),r(r_1),s(s_1);d_1)=\nu(k(k_1),s(s_1))\mu(k(k_1),r(r_1);d_1).
\end{equation}
This yields  a system $(\ref{el2})$ of linear equations in  unknown $\psi(j(j_1),r(r_1),s(s_1);d_1)$  with fixed $r(r_1),s(s_1)$.  Since the matrix $[\nu(k(k_1),j(j_1))]$ of the system $(\ref{el2})$  is a nonsingular,
these equations uniquely determine the $\psi(j(j_1),r(r_1),s(s_1);d_1)$.
$\qed$

\begin{thm}~\label{4t}
The system $(X_d^{d_0},R)$ is an association scheme.
\end{thm}
\proof
We establish the matrix relation as follows
\begin{equation}\label{et1}
   C_rC_s=\sum_{k=0}^d\sum_{\lambda=0}^{h_1}\left\{\sum_{j=0}^{r}\sum_{j'=0}^{h_2}\left[\sum_{r'=0}^{h_3}\left(\psi(j(j'),r(r'),(d-k)(\lambda);d_0)
\sum_{s'=0}^{h_4}\pi(j(j'),r(r'),s(s');d_0)\right)\right]\right\}A_{k}.
\end{equation}
Here $h_1=\min\{d-k,d_0\},h_2=\min\{j,d_0\},h_3=\min\{r,d_0\},h_4=\min\{s,d_0\}.$

In fact, By the definition, for any $D,F\in X_d^{d_0}$, the $(D,F)$-entry of $C_rC_s$ is the number of triples $(A,B,C)\in \widetilde{X_r}\times \widetilde{X_s}\times \widetilde{X_d}$ such that $A\leq D\wedge C$ and $B\leq F\wedge C$ hold.
Let $D\wedge F\in X_{d-k}^{\lambda}$. We assume that $A\wedge F\in X_j^{j'}$ with fixed $j$ and $j'$, by Lemma~\ref{4l-2}, there are $\sum_{r'=0}^{\min\{r,d_0\}}\psi(j(j'),r(r'),(d-k)(\lambda);d_0)$ choices for $A\leq D$. For each $A$, by Lemma~\ref{l3-3}, there are $\sum_{s'=0}^{\min\{s,d_0\}}\pi(j(j'),r(r'),s(s');d_0)$ choices for pairs $(B,C)$ with $B\leq F\wedge C, A\leq C$. Hence, (\ref{et1}) holds.
By Lemma~\ref{4l-1}, $\mathfrak{A}$ is a Bose-Mesner algebra. It implies the assertion.
 $\qed$
\section{$M$-clique}

In this section, we always assume that $l=0$, $X=X_d^0$. By Theorem~\ref{4t}, we obtain a
 association scheme $\mathfrak{X}=(X,\{R_i\}_{0\leq i\leq d})$ in $\mathbb{F}_q^n$.
Next, we will discuss the bound of $M$-clique on association scheme $(X,\{R_i\}_{0\leq i\leq d})$.

 We associate the following linear programming problem, in the real variables $\xi_k(k\in K=\{k\in [0,d]|m_k>1\}):$(see \cite{D23})
 \begin{equation}\label{system}
\left\{
\begin{array}{l}
{\rm min}(f')=1+\sum_{k\in K}(1-b_k)\xi_k,\\
\hbox{subject to}\  \xi_k\geq 0,\  \hbox{all}\  k\in K,\\
\hbox{and to}\ \sum_{k\in K}(P_i(k)-b_kv_i)\xi_k\leq v_i,\ \hbox{all}\ i\in M\backslash \{0\}.
\end{array}
\right.
\end{equation}
Here $P_i(k)$ is eigenvalue of a association scheme, $v_i$ is degree of a association scheme, $b_k=(v_s^{-1}P_s(k))^2,k\in [0,d]$.
\begin{prop}~\label{5p}(\cite{D23})
The problem~(\ref{system}) admits a solution $\min(f')$, and any unicoloured $M$-clique $\alpha$ satisfies
 $\hat{c}_0(\alpha)\leq \min(f')$.
\end{prop}
 By Proposition~\ref{5p}, we obtain the following theorem.
 \begin{thm}
Let $\alpha\in \mathbb{R}^{|X|}$ be a unicoloured $M$-clique on association scheme $\mathfrak{X}=(X,\{R_i\}_{0\leq i\leq d})$ with $M=\{0,1\}$ and
colour $s=1$.  then
$$\hat{c}_0(\alpha)=<\alpha,\alpha>^{-1}<\alpha,1>^2\leq \frac{q^{d+1+n}-3q^{2d+1}+5q^{d+1}-2q^d+2q^{n+d}-3q^{n+1}}{q^{n+d}-q^{2d+1}+2q^{d+1}-q^d-q^{n+1}}.$$
\end{thm}
\proof
   By the (\cite{DD}), the eigenvalues of an association scheme of $\mathfrak{X}=(X,\{R_i\}_{0\leq i\leq d})$ are given by the formula
$$P_i(r)=\sum_{j=0}^i(-1)^{i-j}\left[d-j\atop i-j\right]_q\left[d-r\atop j\right]_q\left[n-d+j-r\atop j\right]_qq^{rj+{i-j\choose 2}}.$$
Pick $\xi_k=\xi\delta_{k,1}$.
It is known that the degree of the $\mathfrak{X}=(X,\{R_i\}_{0\leq i\leq d})$ are given by the formula
$$v_i=q^{i^2}\left[n-d\atop i\right]_q\left[d\atop i\right]_q.$$
By the above formulas, we obtain
$$v_1=\frac{q(q^{n-d}-1)(q^d-1)}{(q-1)^2},$$
$$P_1(1)=-\frac{q^d-1}{q-1}+\frac{q(q^{d-1}-1)(q^{n-d}-1)}{(q-1)^2},$$
For $s=1,i=1$, we have $(P_1(1)-b_1v_1)\xi\leq v_1$, ¼´$(P_1(1)-(v_1^{-1}P_1(1))^2v_1)\xi\leq v_1$. It implies that $(v_1P_1(1)-P_1(1)^2)\xi\leq v_1^2$.
Since
$$v_1-P_1(1)=\frac{q^n-1}{q-1}>0.$$
We have
$$\xi\leq \frac{v_1^2}{p_1(1)(v_1-p_1(1)}.$$
and
$$1-b_1=1-\frac{p_1(1)^2}{v_1^2}=\frac{v_1^2-p_1(1)^2}{v_1^2}>0.$$
It follows that
$$\begin{array}{ll}
\min(f')&=1+(1-b_1)\xi\\
&\leq 1+\frac{v_1+p_1(1)}{p_1(1)}\\
&=\frac{q^{d+1+n}-3q^{2d+1}+5q^{d+1}-2q^d+2q^{n+d}-3q^{n+1}}{q^{n+d}-q^{2d+1}+2q^{d+1}-q^d-q^{n+1}}.
\end{array}$$

By Proposition~\ref{5p}, we have
$$\hat{c}_0(\alpha)\leq  \frac{q^{d+1+n}-3q^{2d+1}+5q^{d+1}-2q^d+2q^{n+d}-3q^{n+1}}{q^{n+d}-q^{2d+1}+2q^{d+1}-q^d-q^{n+1}}.$$
$\qed$

\section*{Acknowledgement}

This research is supported by  NSF of Hebei Province (A2013408009), NSF  of Hebei Education
Department(ZH2012082), the Specialized Research Fund for the Doctoral Program of Higher Education of China (No.20121303110005)
 and the foundation of Langfang Teachers' College
(LSBS201205).


\end{document}